\documentclass[11pt]{article}
\usepackage{amssymb}
\usepackage{amsmath}
\usepackage{amsfonts}
\input amssym.def

\title{\bf Global well-posedness for dissipative Korteweg-de Vries equations}
\author{St\'ephane Vento \\ Universit\'e de Marne-La-Vall\'ee, \\Laboratoire d'Analyse
et de Math\'ematiques Appliqu\'ees,\\ 5 bd. Descartes, Cit\'e
Descartes, Champs-Sur-Marne,\\ 77454 Marne-La-Vall\'ee Cedex 2,
France}

\date{E-mail:\, stephane.vento@univ-mlv.fr}

\numberwithin{equation}{section}

\newtheorem{theorem}{Theorem}[section]
\newtheorem{lemma}{Lemma}[section]
\newtheorem{proposition}{Proposition}[section]

\newtheorem{remark}{Remark}[section]

\def\R{\mathbb{R}}
\def\T{\mathbb{T}}
\def\Z{\mathbb{Z}}
\def\C{\mathcal{C}}
\def\S{\mathcal{S}}

\def\F{\mathcal{F}}
\def\eps{\varepsilon}
\def\supp{\mathop{\rm supp}\nolimits}

\begin{document}
\maketitle {\bf Abstract.}\, This paper is devoted to the
well-posedness for dissipative KdV equations
$u_t+u_{xxx}+|D_x|^{2\alpha}u+uu_x=0$, $0<\alpha\leq 1$. An
optimal bilinear estimate is obtained in Bourgain's type spaces,
which provides global well-posedness in $H^s(\R)$, $s>-3/4$ for
$\alpha\leq1/2$ and $s>-3/(5-2\alpha)$ for $\alpha>1/2$.

\section{Introduction}
We study the initial value problem (IVP) for the dissipative KdV
equations \begin{equation}\label{eq}\left\{\begin{array}{ll}u_t+u_{xxx}+|D_x|^{2\alpha}u+uu_x=0, & t\in\R_+, x\in\R,\\
u(0,x)=u_0(x), & x\in\R,\end{array}\right.\end{equation} with $0<
\alpha\leq 1$ and where $|D_x|^{2\alpha}$ denotes the Fourier
multiplier with symbol $|\xi|^{2\alpha}$. These equations can be
viewed as a combinaison of the KdV equation
\begin{equation}\label{kdv}u_t+u_{xxx}+uu_x=0\end{equation} and Burgers equation
\begin{equation}\label{bur}u_t-u_{xx}+uu_x=0,\end{equation} involving both nonlinear dispersion and
dissipation effects.

The Cauchy problem for the KdV equation has been studied by many
authors. In \cite{MR1215780}, Bourgain introduced new functional
spaces adapted to the linear symbol $\tau-\xi^3$ and showed that
the IVP associated to (\ref{kdv}) is locally well-posed in
$L^2(\R)$. Due to the second conservation law, this result extends
globally in time. Then, working in these spaces, Kenig, Ponce and
Vega obtained local well-posedness in Sobolev spaces $H^s(\R)$ for
$s>-5/8$ in \cite{MR1230283} and for $s>-3/4$ in \cite{MR1329387}.
More recently, global well-posedness was obtained for $s>-3/4$ in
\cite{MR1969209}. It is worth noticing that the index $-3/4$ is
far away from the index $-3/2$ suggested by standard scaling
argument. However, $-3/4$ is indeed the critical index for
well-posedness. In fact, the solution map $u_0\mapsto u$ fails to
be $\mathcal{C}^3$ in $H^s(\R)$, $s<-3/4$ (see \cite{MR1466164})
and $\mathcal{C}^2$ in homogeneous spaces $\dot{H}^s(\R)$,
$s<-3/4$ (see \cite{MR1735881}). Moreover, the bilinear estimate
in $X^{b,s}$ spaces used in \cite{MR1329387} to prove local
well-posedness is sharp with respect of $s$ (see also
\cite{MR1944182}).

Concerning the Cauchy problem for the Burgers equation, the
situation is quite different. By using the strong smoothing effect
of the semigroup related to the heat equation, one can solve
(\ref{bur}) in the Sobolev space given by an heuristic scaling
argument. In \cite{MR1382829}, Dix proved local well-posedness of
(\ref{bur}) in $H^s(\R)$ for $s>-1/2$. Then, this result was
extended to the case $s=-1/2$ in \cite{MR1409926}. Below this
critical index, it has been showed in \cite{MR1382829} that
uniqueness fails.

\vskip 0.5cm

When $\alpha=1/4$, equation (\ref{eq}) models the evolution of the
free surface for shallow water waves damped by viscosity, see
\cite{ott:1432}. When $\alpha=1$, (\ref{eq}) is the so-called
KdV-Burgers equation which models the propagation of weakly
nonlinear dispersive long waves in some contexts when dissipative
effects occur (see \cite{ott:1432}). In \cite{MR1889080}, Molinet
and Ribaud treat the KdV-B equation by working in the usual
Bourgain space related to the KdV equation, considering only the
dispersive part of the equation. They were able to prove global
well-posedness for KdV-B in $H^s(\R)$, $s>-3/4-1/24$, getting a
lower index than the critical indexes for (\ref{kdv}) and
(\ref{bur}). Then, the same authors improved this result in
\cite{MR1918236} by going down to $H^s(\R)$, $s>-1$. The main new
ingredient is the introduction of a new Bourgain space containing
both dispersive and dissipative parts of the equation. For $s<-1$,
the problem is ill-behaved in the sense that the flow map
$u_0\mapsto u$ is not $\mathcal{C}^2$ in $H^s(\R)$.

Concerning the case $0<\alpha<1$ in (\ref{eq}), Molinet and Ribaud
established in \cite{MR1889080} the global well-posedness for data
in $H^s(\R)$, $s>-3/4$, whatever the value of $\alpha$. On the
other hand, ill-posedness is known for (\ref{eq}) in
$\dot{H}^s(\R)$, $s<(\alpha-3)/(2(2-\alpha))$, see
\cite{MR1918236}.

\vskip 0.5cm

In this paper we improve the results obtained in \cite{MR1889080}.
We show that the Cauchy problem (\ref{eq}) is globally well-posed
in $H^s(\R)$, $s>s_\alpha$ with
$$s_\alpha=\left\{\begin{array}{ll}-3/4 & \mathrm{ if
}\ 0<\alpha\leq 1/2,\\ -\frac{3}{5-2\alpha} & \mathrm{ if }\
1/2<\alpha\leq 1\end{array}\right..$$ Of course the case
$\alpha\leq 1/2$ is well-known, but our general proofs contain
this result. In suitable $X^{b,s}$ spaces, we are going to perform
a fixed point argument on the integral formulation of (\ref{eq}).
This will be achieved by deriving a bilinear estimate in these
spaces. By Plancherel's theorem and duality, it reduces to
estimating a weighted convolution of $L^2$ functions. To recover
the lost derivative in the nonlinear term $\partial_x(u^2)$, we
take advantage of the well-known algebraic smoothing relation
(\ref{smooth}) combined with several methods. On one hand, we use
Strichartz's type estimates and some techniques introduced in
\cite{MR1329387}. On the other hand, these techniques are not
sufficient in certain regions to go down below $-3/4$ and we are
lead tu use dyadic decomposition. In \cite{MR1854113}, Tao studied
some nonlinear dispersive equations like KdV, Schr\"{o}dinger or
wave equation by using such dyadic decomposition and
orthogonality. He obtained sharp estimates on dyadic blocs, which
leads to multilinear estimates in the $X^{b,s}$ spaces, usable in
many contexts. Note that very recently, such a method was
exploited in \cite{chen-2007} for the dissipative modified-KdV
equation $u_t+u_{xxx}+|D_x|^{2\alpha}u+u^2u_x=0$.

It is worth pointing out that (\ref{eq}) has no scaling
invariance, so it is not clear that $s_\alpha$ is the critical
index for well-posedness. However, the fundamental bilinear
estimate used to prove our results is optimal.

\subsection{Notations}
For two positive reals $A$ and $B$, we write $A\lesssim B$ if
there exists a constant $C>0$ such that $A\leq CB$. When this
constant is supposed to be sufficiently small, we write $A\ll B$.
Similarly, we use the notations $A\gtrsim B$, $A\sim B$ and $A\gg
B$. When $x\in\R$, $x^+$ denotes its positive part $\max(0,x)$.
For $f\in \S'(\R^N)$, we define its Fourier transform $\F(f)$
 (or $\widehat{f}$) by
$$\F f(\xi)=\int_{\R^N}e^{-i\langle x,\xi\rangle}f(x)dx.$$
The Lebesgue spaces are endowed with the norm
$$\|f\|_{L^p(\R^N)}=\Big(\int_{\R^N}|f(x)|^pdx\Big)^{1/p},\quad
1\leq p<\infty$$ with the usual modification for $p=\infty$. We
also consider the space-time Lebesgue spaces $L^p_xL^q_t$ defined
by
$$\|f\|_{L^p_xL^q_t}=\Big\|\|f\|_{L^p_x(\R)}\Big\|_{L^q_t(\R)}.$$
For $b,s\in\R$, we define the Sobolev spaces $H^s(\R)$ and their
space-time versions $H^{b,s}(\R^2)$ by the norms
$$\|f\|_{H^s}=\Big(\int_\R\langle\xi\rangle^{2s}|\widehat{f}(\xi)|^2d\xi\Big)^{1/2},$$
$$\|u\|_{H^{b,s}}=\Big(\int_{\R^2}\langle\tau\rangle^{2b}\langle\xi\rangle^{2s}|\widehat{u}(\tau,\xi)|^2d\tau
d\xi\Big)^{1/2},$$ with $\langle\cdot\rangle=(1+|\cdot|^2)^{1/2}$.

Let $U(t)\varphi$ denote the solution of the Airy equation
$$\left\{\begin{array}{ll}u_t+u_{xxx}=0\\u(0)=\varphi\end{array}\right.,$$
that is, $$\forall t\in\R,\
\F_x(U(t)\varphi)(\xi)=\exp(i\xi^3t)\widehat{\varphi}(\xi),\quad
\varphi\in\S'.$$ In \cite{MR1918236}, Molinet and Ribaud
introduced the function spaces $X^{b,s}_\alpha$ related to the
linear symbol $i(\tau-\xi^3)+|\xi|^{2\alpha}$ and defined by the
norm
$$\|u\|_{X^{b,s}_\alpha}=\big\|\langle
i(\tau-\xi^3)+|\xi|^{2\alpha}\rangle^b\langle\xi\rangle^s\widehat{u}\big\|_{L^2(\R^2)}.$$
Note that since
$\F(U(-t)u)(\tau,\xi)=\widehat{u}(\tau+\xi^3,\xi)$, we can
re-express the norm of $X^{b,s}_\alpha$ as \begin{eqnarray*}
\|u\|_{X^{b,s}_\alpha} &=& \big\|\langle
i\tau+|\xi|^{2\alpha}\rangle^b\langle\xi\rangle^s\widehat{u}(\tau+\xi^3,\xi)\big\|_{L^2(\R^2)}\\
&=& \big\|\langle
i\tau+|\xi|^{2\alpha}\rangle^b\langle\xi\rangle^s\F(U(-t)u)(\tau,\xi)\big\|_{L^2(\R^2)}\\
&\sim & \|U(-t)u\|_{H^{b,s}}+\|u\|_{L^2_tH^{s+2\alpha b}_x}.
\end{eqnarray*}
We will also work in the restricted spaces $X^{b,s}_{\alpha,T}$,
$T\geq 0$, equipped with the norm
$$\|u\|_{X^{b,s}_{\alpha,T}}=\inf_{w\in
X^{b,s}_\alpha}\{\|w\|_{X^{b,s}_\alpha},\ w(t)=u(t)\textrm{ on }
[0,T]\}.$$ Finally, we denote by $W_\alpha$ the semigroup
associated with the free evolution of (\ref{eq}), $$\forall t\geq
0,\
\F_x(W_\alpha(t)\varphi)(\xi)=\exp[-|\xi|^{2\alpha}t+i\xi^3t]\widehat{\varphi}(\xi),\
\varphi\in\S',$$ and we extend $W_\alpha$ to a linear operator
defined on the whole real axis by setting
\begin{equation}\label{walpha}\forall t\in\R,\
\F_x(W_\alpha(t)\varphi)(\xi)=\exp[-|\xi|^{2\alpha}|t|+i\xi^3t]\widehat{\varphi}(\xi),\
\varphi\in\S'.\end{equation}

\subsection{Main results}
Let us first state our crucial bilinear estimate.
\begin{theorem}\label{thbil} Given $s>s_\alpha$, there exist $\nu,\delta>0$
such that for any  $u,v\in X^{1/2,s}_\alpha$ with compact support
in $[-T,+T]$,
\begin{equation}\label{bil}\|\partial_x(uv)\|_{X^{-1/2+\delta,s}_\alpha}\lesssim
T^\nu
\|u\|_{X^{1/2,s}_\alpha}\|v\|_{X^{1/2,s}_\alpha}.\end{equation}
\end{theorem}\label{opt}
This result is optimal in the following sense.
\begin{theorem}For all $s\leq s_\alpha$ and $\nu,\delta>0$,  there exist $u,v\in X^{1/2,s}_\alpha$ with compact support in
$[-T,+T]$ such that the estimate (\ref{bil}) fails.
\end{theorem}

Let $\psi$ be a cutoff function such that
$$\psi\in\mathcal{C}^\infty_0(\R),\quad \supp \psi\subset
[-2,2],\quad \psi\equiv 1\textrm{ on }[-1,1],$$ and define
$\psi_T(\cdot)=\psi(\cdot/T)$ for all $T>0$. By Duhamel's
principle, the solution to the problem (\ref{eq}) can be locally
written in the integral form as
\begin{equation}\label{int}u(t)=\psi(t)\Big[W_\alpha(t)u_0-\frac{\chi_{\R_+}(t)}{2}\int_0^tW_\alpha(t-t')\partial_x(\psi_T^2(t')
u^2(t'))dt'\Big].
\end{equation}
Clearly, if $u$ is a solution of (\ref{int}) on $[-T,+T]$, then
$u$ solves (\ref{eq}) on $[0,T/2]$.

 As a consequence of Theorem
\ref{thbil} together with linear estimates of Section
\ref{subsec-lin}, we obtain the following global well-posedness
result.

\begin{theorem}\label{main} Let $\alpha\in (0,1]$ and $u_0\in H^s(\R)$ with
$s>s_\alpha$. Then for any $T>0$, there exists a unique solution
$u$ of (\ref{eq}) in $$Z_T=\mathcal{C}([0,T], H^s(\R))\cap
X^{1/2,s}_{\alpha,T}.$$ Moreover, the map $u_0\mapsto u$ is smooth
from $H^s(\R)$ to $Z_T$ and $u$ belongs to $\mathcal{C}((0,T],
H^\infty(\R))$.
\end{theorem}

\begin{remark}
Actually, we shall prove Theorems \ref{thbil} and \ref{main} in
the most difficult case. In the sequel we assume
\begin{equation}\label{conds}
\left\{\begin{array}{ll}s_\alpha<s<-1/2 & \textrm{ if }\ \alpha\leq 1/2,\\
s_\alpha<s<-3/4 &\textrm{ if }\ \alpha>1/2.\end{array}\right.
\end{equation}
\end{remark}

\begin{remark}
Theorem \ref{main} is known to be sharp in the case $\alpha=0$
(KdV equation) and in the case $\alpha=1$ (KdV-B equation). On the
other hand, as far as we know, most of nonlinear equations for
which the multilinear estimate fails in the related $X^{b,s}$
space
 are ill-posed in $H^s$. Therefore it is reasonable to
conjecture that $s_\alpha$ is really the critical index for
(\ref{eq}). The fact that $s_\alpha=-3/4$ for $\alpha\leq 1/2$
could mean that the dissipative part in (\ref{eq}), when becoming
small enough, has no effect on the low regularity of the equation.
\end{remark}

\begin{remark}
It is an interesting problem to consider the periodic dissipative
KdV equation
\begin{equation}\label{period}\left\{\begin{array}{ll}u_t+u_{xxx}-|D_x|^{2\alpha}u+uu_x=0, & t\in\R_+, x\in\T,\\
u(0,x)=u_0(x), & x\in\T,\end{array}\right..\end{equation}
Concerning the KdV equation on $\T$, global well-posedness is
known in $H^{1/2}(\T)$ (see \cite{MR1969209}) and the result is
optimal (see \cite{MR1466164}). For KdV-B, it is established in
\cite{MR1918236} that the indexes of the critical spaces are the
same on the real line and on the circle. We believe that working
in the space $\widetilde{X}^{1/2,s}_\alpha$ endowed with the norm
$$\|u\|_{\widetilde{X}^{1/2,s}_\alpha}=\Big(\sum_{n\in\Z}\langle
n\rangle^{2s}\int_\R \langle i(\tau-n^3)+|n|^{2\alpha}\rangle
|\widehat{u}(\tau,n)|^2d\tau\Big)^{1/2},$$ and using Tao's
$[k;Z]$-multiplier norm estimates \cite{MR1854113}, one could get
well-posedness results for the IVP (\ref{period}) in $H^s(\T)$,
$s>\tilde{s}_\alpha$ with $\tilde{s}_0=-1/2$ and $\tilde{s}_1=-1$.
We do not pursue this issue here.
\end{remark}

The remainder of this paper is organized as follows. In Section
\ref{sec-pre}, we recall some linear estimates on the operators
$W_\alpha$ and $L_\alpha$, and we introduce Tao's
$[k;Z]$-multiplier norm estimates. Section \ref{sec-bil} is
devoted to the proof of the bilinear estimate (\ref{bil}). Theorem
\ref{main} is established in Section \ref{sec-main}. Finally, we
show the optimality of (\ref{bil}) in Section \ref{sec-opt}.

\section*{Acknowledgment}
The author would like to express his gratitude to Francis Ribaud
for his availability  and his constant encouragements.

\section{Preliminaries}
\label{sec-pre}
\subsection{Linear estimates}
\label{subsec-lin} In this subsection, we collect together several
linear estimates on the operators $W_\alpha$ introduced in
(\ref{walpha}) and $L_\alpha$ defined by $$L_\alpha: f\mapsto
\chi_{\R_+}(t)\psi(t)\int_0^tW_\alpha(t-t')f(t')dt'.$$ All the
results stated here were proved in \cite{MR1918236} for $\alpha=1$
and in \cite{chen-2007} for the general case.

\begin{lemma} For all $s\in\R$ and all $\varphi\in H^s(\R)$,
\begin{equation}\label{lin-free}\|\psi(t)W_\alpha(t)\varphi\|_{X^{1/2,s}_\alpha}\lesssim
\|\varphi\|_{H^s}.\end{equation}
\end{lemma}

\begin{lemma} Let $s\in\R$.
\begin{enumerate}
\item[(a)] For all $v\in\mathcal{S}(\R^2)$,
\begin{eqnarray*}\lefteqn{\Big\|\chi_{\R_+}(t)\psi(t)\int_0^tW_\alpha(t-t')v(t')dt'\Big\|_{X^{1/2,s}_\alpha}}\\
& \lesssim
\|v\|_{X^{-1/2,s}_\alpha}+\Big(\int_\R\langle\xi\rangle^{2s}\big(\int_\R
\frac{|\widehat{v}(\tau+\xi^3,\xi)|}{\langle
i\tau+|\xi|^{2\alpha}\rangle}d\tau\big)d\xi\Big)^{1/2}.
\end{eqnarray*}
\item[(b)]For all $0<\delta<1/2$ and all $v\in
X^{-1/2+\delta,s}_\alpha$,
\begin{equation}\label{lin-for}\Big\|\chi_{\R_+}(t)\psi(t)\int_0^tW_\alpha(t-t')v(t')dt'\Big\|_{X^{1/2,s}_\alpha}\lesssim
\|v\|_{X^{-1/2+\delta,s}_\alpha}.\end{equation}
\end{enumerate}
\end{lemma}
To globalize our solution, we will need the next lemma.
\begin{lemma}\label{glob} Let $s\in\R$ and $\delta>0$. Then for any
$f\in X^{-1/2+\delta,s}_\alpha$, $$t\longmapsto
\int_0^tW_\alpha(t-t')f(t')dt'\in\mathcal{C}(\R_+,H^{s+2\alpha\delta}).$$
Moreover, if $(f_n)$ is a sequence satisfying $f_n\rightarrow 0$
in  $X^{-1/2+\delta,s}_\alpha$, then
$$\Big\|\int_0^tW_\alpha(t-t')f_n(t')dt'\Big\|_{L^\infty(\R_+,H^{s+2\alpha\delta})}\longrightarrow
0.$$
\end{lemma}

Finally, we recall the following $L^4$ Strichartz's type estimate
showed in \cite{MR1230283,MR1889080}.

\begin{lemma}\label{stri} Let $f\in L^2(\R^2)$ with compact support (in time) in $[-T,+T]$. For $0\leq \theta\leq 1/8$
and $\rho>3/8$, there exists $\nu>0$ such that
$$\Big\|\mathcal{F}^{-1}\Big(\frac{\langle\xi\rangle^\theta
\widehat{f}(\tau,\xi)}{\langle\tau-\xi^3\rangle^\rho}\Big)\Big\|_{L^4_{xt}}\lesssim
T^\nu \|f\|_{L^2_{xt}}.$$
\end{lemma}

\subsection{Tao's $[k;Z]$-multipliers}
\label{subsec-tao} Now we turn to Tao's $[k;Z]$-multiplier norm
estimates. For more details, please refer to \cite{MR1854113}.

Let $Z$ be any abelian additive group with an invariant measure
$d\xi$. For any integer $k\geq 2$ we define the hyperplane
$$\Gamma_k(Z)=\{(\xi_1,...,\xi_k)\in Z^k : \xi_1+...+\xi_k=0\}$$
which is endowed with the measure
$$\int_{\Gamma_k(Z)}f=\int_{Z^{k-1}}f(\xi_1,...,\xi_{k-1},-(\xi_1+...+\xi_{k-1}))d\xi_1...d\xi_{k-1}.$$
A $[k;Z]$-multiplier is defined to be any function
$m:\Gamma_k(Z)\rightarrow \mathbb{C}$. The multiplier norm
$\|m\|_{[k;Z]}$ is defined to be the best constant such that the
inequality
\begin{equation}\label{def-norm-mult}\Big|\int_{\Gamma_k(Z)}m(\xi)\prod_{j=1}^kf_j(\xi_j)\Big|\leq
\|m\|_{[k;Z]}\prod_{j=1}^k\|f_j\|_{L^2(Z)}\end{equation} holds for
all test functions $f_1,...,f_k$ on $Z$. In other words,
$$\|m\|_{[k;Z]}=\sup_{\substack{f_j\in\S(Z)\\\|f_j\|_{L^2(Z)}\leq
1}}\Big|\int_{\Gamma_k(Z)}m(\xi)\prod_{j=1}^kf_j(\xi_j)\Big|.$$ In
his paper \cite{MR1854113}, Tao used the following notations.
Capitalized variables $N_j$, $L_j$ ($j=1,...,k$) are presumed to
be dyadic, i.e. range over numbers of the form $2^\ell$,
$\ell\in\Z$. In this paper, we only consider the case $k=3$, which
corresponds to the quadratic nonlinearity in the equation. It will
be convenient to define the quantities $N_{max}\geq N_{med}\geq
N_{min}$ to be the maximum, median and minimum of $N_1,N_2,N_3$
respectively. Similarly, define $L_{max}\geq L_{med}\geq L_{min}$
whenever $L_1,L_2,L_3>0$. The quantities $N_j$ will measure the
magnitude of frequencies of our waves, while $L_j$ measures how
closely our waves approximate a free solution.

Here we consider $[3,\R\times\R]$-multipliers and we parameterize
$\R\times\R$ by $(\tau,\xi)$ endowed with the Lebesgue measure
$d\tau d\xi$. If $\tau, \tau_1, \xi,\xi_1$ are given, we set
\begin{equation}
\sigma=\sigma(\tau,\xi)=\tau-\xi^3,\quad
\sigma_1=\sigma(\tau_1,\xi_1),\quad
\sigma_2=\sigma(\tau-\tau_1,\xi-\xi_1).
\end{equation}
From the identity $\sigma_1+\sigma_2-\sigma=3\xi\xi_1(\xi-\xi_1)$
one can deduce the well-known smoothing relation
\begin{equation}\label{smooth}\max(|\sigma|,|\sigma_1|,|\sigma_2|)\geq
|\xi\xi_1(\xi-\xi_1)|\end{equation} which will be extensively used
in Section \ref{sec-bil}.

 By a dyadic
decomposition of the variables $\xi_1$, $\xi_2=-\xi$,
$\xi_3=\xi-\xi_1$, and $\sigma_1$, $\sigma_2$, $\sigma_3=-\sigma$,
 we are lead to consider
\begin{equation}\label{mult}\Big\|\prod_{j=1}^3\chi_{|\xi_j|\sim N_j}\chi_{|\sigma_j|\sim
L_j}\Big\|_{[3;\R\times\R]}.\end{equation} We can now state the
fundamental dyadic estimates for the KdV equation on the real line
(\cite{MR1854113}, Proposition 6.1).

\begin{lemma}\label{tao} Let $N_1, N_2, N_3, L_1, L_2, L_3$ satisfying
$$N_{max}\sim N_{med},$$ $$L_{max}\sim \max(N_1N_2N_3, L_{med}).$$
\begin{enumerate}
\item ((++) Coherence) If $N_{max}\sim N_{min}$ and $L_{max}\sim
N_1N_2N_3$ then we have \begin{equation}(\ref{mult})\lesssim
L_{min}^{1/2}N_{max}^{-1/4}L_{med}^{1/4}.\end{equation} \item
((+-) Coherence) If $N_2\sim N_3\gg N_1$ and $N_1N_2N_3\sim
L_1\gtrsim L_2, L_3$ then
\begin{equation}\label{+-} (\ref{mult})\lesssim
L_{min}^{1/2}N_{max}^{-1}\min\big(N_1N_2N_3,
\frac{N_{max}}{N_{min}}L_{med}\big)^{1/2}.
\end{equation}
Similarly for permutations. \item In all other cases, we have
\begin{equation}(\ref{mult})\lesssim
L_{min}^{1/2}N_{max}^{-1}\min(N_1N_2N_3,
L_{med})^{1/2}.\end{equation}
\end{enumerate}
\end{lemma}

Because only one region needs to be controlled by using a dyadic
approach, we just require the (+-) coherence case. On the other
hand, these estimates are sharp. In particular, testing
(\ref{def-norm-mult}) with
\begin{eqnarray*} f_1(\tau,\xi) &=& \chi_{|\xi|\sim N_1;
|\tau-3N_2^2\xi|\lesssim N_1^2N_2},\\ f_2(\tau,\xi) &=&
\chi_{|\xi-N_2|\lesssim N_1; |\sigma|\lesssim L_2},\\
f_3(\tau,\xi) &=& \chi_{|\xi+N_2|\lesssim N_1; |\sigma|\lesssim
L_3},
\end{eqnarray*}
one obtain the optimality of bound (\ref{+-}) in the case $N_2\sim
N_3\gtrsim N_1$ and $N_1N_2N_3\sim L_1\gtrsim L_2\gtrsim L_3$.
This will be crucial in the proof of Theorem \ref{opt}.

\section{Bilinear estimate}
\label{sec-bil} In this section, we derive the bilinear estimate
(\ref{bil}). To get the required contraction factor $T^\nu$ in our
estimates, the next lemma is very useful (see \cite{MR1491547}).

\begin{lemma}\label{contract} Let $f\in L^2(\R^2)$ with compact
support (in time) in $[-T, +T]$. For any $\theta>0$, there exists
$\nu=\nu(\theta)>0$ such that
$$\Big\|\mathcal{F}^{-1}\Big(\frac{\widehat{f}(\tau,\xi)}{\langle\tau-\xi^3\rangle^\theta}\Big)\Big\|_{L^2_{xt}}\lesssim
T^\nu\|f\|_{L^2_{xt}}.$$
\end{lemma}
We will also need the following elementary calculus inequalities.
\begin{lemma}\label{calc}
\begin{enumerate}
\item[(a)] For $b,b'\in]\frac 14,\frac 12[$ and
$\alpha,\beta\in\R$,
\begin{equation}\label{calc1}\int_\R\frac{dx}{\langle
x-\alpha\rangle^{2b}\langle x-\beta\rangle^{2b'}}\lesssim
\frac{1}{\langle \alpha-\beta\rangle^{2b+2b'-1}}.\end{equation}
\item[(b)] For $b,b'\in]\frac 14,\frac 12[$ and
$\alpha,\beta\in\R$, \begin{equation}\label{calc2}\int_{|x|\leq
|\beta|}\frac{dx}{\langle
x\rangle^{2b+2b'-1}\sqrt{|\alpha-x|}}\lesssim
\frac{\langle\beta\rangle^{2(1-b-b')}}{\langle\alpha\rangle^{1/2}}.\end{equation}
\end{enumerate}
\end{lemma}
\vspace{1cm} {\bf Proof of Theorem \ref{thbil} :} By duality,
(\ref{bil}) is equivalent to
$$\Big|\int_{\R^2}\partial_x(uv)w\Big|\lesssim
T^\nu\|w\|_{X^{1/2-\delta,-s}_\alpha}\|u\|_{X^{1/2,s}_\alpha}\|v\|_{X^{1/2,s}_\alpha}$$
for all $w\in X^{1/2,s}_\alpha$, and setting
\begin{eqnarray*} \widehat{f}(\tau,\xi) &=& \langle
i(\tau-\xi^3)+|\xi|^{2\alpha}\rangle^{1/2}\langle\xi\rangle^s\widehat{u}(\tau,\xi),\\
\widehat{g}(\tau,\xi) &=& \langle
i(\tau-\xi^3)+|\xi|^{2\alpha}\rangle^{1/2}\langle\xi\rangle^s\widehat{v}(\tau,\xi),\\
\widehat{h}(\tau,\xi) &=& \langle
i(\tau-\xi^3)+|\xi|^{2\alpha}\rangle^{1/2-\delta}\langle\xi\rangle^{-s}\widehat{w}(\tau,\xi),
\end{eqnarray*}
it is equivalent to show that
$$I=\int_{\R^4}K(\tau,\tau_1,\xi,\xi_1)\widehat{h}(\tau,\xi)\widehat{f}(\tau_1,\xi_1)\widehat{g}(\tau-\tau_1,\xi-\xi_1)d\tau
d\tau_1 d\xi d\xi_1\lesssim
T^\nu\|f\|_{L^2_{xt}}\|g\|_{L^2_{xt}}\|h\|_{L^2_{xt}}$$ with
$$K=\frac{|\xi|\langle\xi\rangle^{s}}{\langle
i\sigma+|\xi|^{2\alpha}\rangle^{1/2-\delta}}
\frac{\langle\xi_1\rangle^{-s}}{\langle
i\sigma_1+|\xi_1|^{2\alpha}\rangle^{1/2}}\frac{\langle\xi-\xi_1\rangle^{-s}}{\langle
i\sigma_2+|\xi-\xi_1|^{2\alpha}\rangle^{1/2}}.$$ By Fubini's
theorem, we can always assume $\widehat{f}, \widehat{g},
\widehat{h}\geq 0$. By symmetry, one can reduce the integration
domain of $I$ to
$\Omega=\{(\tau,\tau_1,\xi,\xi_1)\in\R^4,|\sigma_1|\geq
|\sigma_2|\}$. Split $\Omega$ into four regions,
\begin{eqnarray*}\Omega_1 &=&
\{(\tau,\tau_1,\xi,\xi_1)\in\Omega : |\xi_1|\leq 1,|\xi-\xi_1|\leq
1\},\\
\Omega_2 &=& \{(\tau,\tau_1,\xi,\xi_1)\in\Omega : |\xi_1|\leq
1,|\xi-\xi_1|\geq 1\},\\ \Omega_3 &=&
\{(\tau,\tau_1,\xi,\xi_1)\in\Omega : |\xi_1|\geq 1,|\xi-\xi_1|\leq
1\},\\ \Omega_4 &=& \{(\tau,\tau_1,\xi,\xi_1)\in\Omega :
|\xi_1|\geq 1,|\xi-\xi_1|\geq 1\}.
\end{eqnarray*}

\textbf{Estimate in $\Omega_1$}\\
Using Cauchy-Schwarz inequality and Lemma \ref{contract}, we
easily obtain
\begin{equation}\label{est}\begin{split}I_1 &\lesssim
\sup_{\tau,\xi}\Big[\frac{|\xi|\langle\xi\rangle^{s}}{\langle
i\sigma+|\xi|^{2\alpha}\rangle^{1/2-\delta/2}}\Big(\int_{\widetilde{\Omega}_1}\frac{\langle\xi_1\rangle^{-2s}
\langle\xi-\xi_1\rangle^{-2s}}{\langle
i\sigma_1+|\xi_1|^{2\alpha}\rangle\langle
i\sigma_2+|\xi-\xi_1|^{2\alpha}\rangle}d\tau_1
d\xi_1\Big)^{1/2}\Big]\\ &\times
T^\nu\|f\|_{L^2_{xt}}\|g\|_{L^2_{xt}}\|h\|_{L^2_{xt}}\end{split}\end{equation}
with $\widetilde{\Omega}_1=\{(\tau_1,\xi_1) : \exists
\tau,\xi\in\R, (\tau,\tau_1,\xi,\xi_1)\in\Omega_1\}$.\\
In $\Omega_1$, one has $|\xi|\leq 2$ ans thus if $K_1$ denotes the
term between brackets in (\ref{est}),
\begin{eqnarray*}K_1 &\lesssim &
\Big(\int_{\widetilde{\Omega}_1}\frac{d\tau_1 d\xi_1}{\langle
\sigma_1\rangle\langle\sigma_2\rangle}\Big)^{1/2}\lesssim
\Big(\int_{|\xi_1|\leq
1}\Big(\int_\R\frac{d\tau_1}{\langle\sigma_2\rangle^{2}}\Big)d\xi_1\Big)^{1/2}\lesssim
1.
\end{eqnarray*}

\textbf{Estimate in $\Omega_2$}\\
We split $\Omega_2$ into
\begin{eqnarray*}
\Omega_{21}&=&\{(\tau,\tau_1,\xi,\xi_1)\in\Omega_2 :
|\sigma|\geq|\sigma_1|\},\\
\Omega_{22}&=&\{(\tau,\tau_1,\xi,\xi_1)\in\Omega_2 :
|\sigma_1|\geq|\sigma|\}.
\end{eqnarray*}

\underline{Estimate in $\Omega_{21}$ :} Note that in this region,
$|\xi-\xi_1|\sim \langle\xi-\xi_1\rangle$ and
$\langle\xi\rangle\gtrsim \langle\xi-\xi_1\rangle$. Thus using
(\ref{est}) as well as (\ref{smooth}), it follows that
\begin{eqnarray*}
K_{21} & \lesssim &
\langle\xi\rangle^{1/2+s+\delta/2}\Big(\int_{\widetilde{\Omega}_{21}}\frac{\langle\xi-\xi_1\rangle^{-2s-1+\delta}}
{|\xi_1|^{1-\delta}\langle\sigma_1\rangle\langle\sigma_2\rangle^\eps\langle\xi-\xi_1\rangle^{2\alpha(1-\eps)}}d\tau_1
d\xi_1\Big)^{1/2}\\
 &\lesssim & \Big(\int_{|\xi_1|\leq
1}\frac{\langle\xi-\xi_1\rangle^{-2\alpha(1-\eps)+2\delta}}{|\xi_1|^{1-\delta}}\Big(\int_\R\frac{d\tau_1}{\langle\sigma_1\rangle
\langle\sigma_2\rangle^\eps}\Big)d\xi_1\Big)^{1/2} \\ &\lesssim &
\Big(\int_{|\xi_1|\leq
1}\frac{d\xi_1}{|\xi_1|^{1-\delta}}\Big)^{1/2}\\ &\lesssim & 1.
\end{eqnarray*}

\underline{Estimate in $\Omega_{22}$ :} By similar arguments, we
estimate
\begin{eqnarray*}K_{22} &\lesssim &
\frac{|\xi|\langle\xi\rangle^{s}}{\langle\xi\rangle^{\alpha(1-\delta)}}\Big(\int_{\widetilde{\Omega}_{22}}\frac{\langle\xi-\xi_1\rangle^{-2s}}
{|\xi\xi_1(\xi-\xi_1)|^{1-\eps}\langle\sigma_1\rangle^\eps\langle\sigma_2\rangle}d\tau_1d\xi_1\Big)^{1/2}\\
&\lesssim &
\langle\xi\rangle^{1/2+s-\alpha(1-\delta)+\eps/2}\Big(\int_{|\xi_1|\leq
1}\frac{\langle\xi-\xi_1\rangle^{-2s-1+\eps}}{|\xi_1|^{1-\eps}}d\xi_1\Big)^{1/2}\\
&\lesssim & \langle\xi\rangle^{-\alpha(1-\delta)+\eps}\\ &\lesssim
& 1.
\end{eqnarray*}

\textbf{Estimate in $\Omega_3$}\\
By symmetry, the desired bound in this region can be obtained in
the same way.

\textbf{Estimate in $\Omega_4$}\\
Divide $\Omega_4$ into
\begin{eqnarray*}
\Omega_{41}&=&\{(\tau,\tau_1,\xi,\xi_1)\in\Omega_4 :
|\sigma_1|\geq|\sigma|\},\\
\Omega_{42}&=&\{(\tau,\tau_1,\xi,\xi_1)\in\Omega_4 :
|\sigma|\geq|\sigma_1|\}.
\end{eqnarray*}

\textbf{Estimate in $\Omega_{41}$}\\
We write
$\Omega_{41}=\Omega_{411}\cup\Omega_{412}\cup\Omega_{413}$ with
\begin{eqnarray*}\Omega_{411}&=&\{(\tau,\tau_1,\xi,\xi_1)\in\Omega_{41}
: |\xi_1|\leq 100|\xi|\},\\
\Omega_{412}&=&\{(\tau,\tau_1,\xi,\xi_1)\in\Omega_{41} :
|\xi_1|\geq 100|\xi|, 3|\xi\xi_1(\xi-\xi_1)|\leq
\frac 12 |\sigma_1|\},\\
\Omega_{413}&=&\{(\tau,\tau_1,\xi,\xi_1)\in\Omega_{41} :
|\xi_1|\geq 100|\xi|, 3|\xi\xi_1(\xi-\xi_1)|\geq \frac 12
|\sigma_1|\}.
\end{eqnarray*}

\underline{Estimate in $\Omega_{411}$ :} We have
$\langle\xi_1\rangle\lesssim\langle\xi\rangle$ and
$\langle\xi-\xi_1\rangle\lesssim \langle\xi\rangle$ thus with
(\ref{smooth}), we deduce for $0<\lambda<1$,
\begin{eqnarray*}K_{411} &\lesssim &
\frac{\langle\xi\rangle^{1/2+s}}{\langle\sigma\rangle^{\lambda/2}\langle\xi\rangle^{\alpha(1-\lambda-2\delta)}}
\langle\xi_1\rangle^{-s-1/2}
\frac{\langle\xi-\xi_1\rangle^{-s-1/2}}{\langle\sigma_2\rangle^{\lambda/2}\langle\xi-\xi_1\rangle^{\alpha(1-\lambda)}}\\
&\lesssim &
\frac{\langle\xi\rangle^{[-s/2-1/4-\alpha(1-\lambda-\delta)]^+}}{\langle\sigma\rangle^{\lambda/2}}
\frac{\langle\xi-\xi_1\rangle^{[-s/2-1/4-\alpha(1-\lambda-\delta)]^+}}{\langle\sigma_2\rangle^{\lambda/2}}.
\end{eqnarray*}
Consequently, using Plancherel's theorem and H\"{o}lder
inequality,
\begin{eqnarray*}I_{411} &\lesssim &
\int_{\R^2}\mathcal{F}^{-1}\Big(\frac{\langle\xi\rangle^{[-s/2-1/4-\alpha(1-\lambda-\delta)]^+}\widehat{h}}
{\langle\sigma\rangle^{\lambda/2}}\Big)\mathcal{F}^{-1}(\widehat{f})
\mathcal{F}^{-1}\Big(\frac{\langle\xi\rangle^{[-s/2-1/4-\alpha(1-\lambda-\delta)]^+}\widehat{g}}
{\langle\sigma\rangle^{\lambda/2}}\Big)dtdx\\ &\lesssim &
\Big\|\mathcal{F}^{-1}\Big(\frac{\langle\xi\rangle^{[-s/2-1/4-\alpha(1-\lambda-\delta)]^+}\widehat{h}}
{\langle\sigma\rangle^{\lambda/2}}\Big)\Big\|_{L^4_{xt}}\|f\|_{L^2_{xt}}\Big\|
\mathcal{F}^{-1}\Big(\frac{\langle\xi\rangle^{[-s/2-1/4-\alpha(1-\lambda-\delta)]^+}\widehat{g}}
{\langle\sigma\rangle^{\lambda/2}}\Big)\Big\|_{L^4_{xt}}.
\end{eqnarray*}
Now we choose $\lambda=3/4+\delta$ so that $\lambda/2>3/8$ and for
$s>-3/4-\alpha/2$ and $\delta>0$ small enough,
$$-\frac s2-\frac 14-\alpha(1-\lambda-\delta)=\frac 18-\frac 12\big(s+\frac 34+\frac\alpha 2\big)
+2\alpha\delta\leq \frac 18.$$
 Hence with help of Lemma \ref{stri},
$I_{411}\lesssim
T^\nu\|f\|_{L^2_{xt}}\|g\|_{L^2_{xt}}\|h\|_{L^2_{xt}}$.

\vskip 0.5cm \underline{Estimate in $\Omega_{412}$ :} Using the
same arguments that for (\ref{est}), we show that
\begin{equation}\label{est1}\begin{split}I_{412} &\lesssim
\sup_{\tau_1,\xi_1}\Big[\frac{\langle\xi_1\rangle^{-s}}{\langle
i\sigma_1+|\xi_1|^{2\alpha}\rangle^{1/2}}\Big(\int_{\widetilde{\Omega}_{412}}\frac{|\xi|^2\langle\xi\rangle^{2s}
\langle\xi-\xi_1\rangle^{-2s}}{\langle
i\sigma+|\xi|^{2\alpha}\rangle^{1-\delta}\langle
i\sigma_2+|\xi-\xi_1|^{2\alpha}\rangle}d\tau
d\xi\Big)^{1/2}\Big]\\ &\times
T^\nu\|f\|_{L^2_{xt}}\|g\|_{L^2_{xt}}\|h\|_{L^2_{xt}}\end{split}\end{equation}
with $\widetilde{\Omega}_{412}=\{(\tau,\xi)
:\exists\tau_1,\xi_1\in\R,
(\tau,\tau_1,\xi,\xi_1)\in\Omega_{422}\}$. Moreover, we easily
check that in $\Omega_{412}$,
$$|\sigma_1+3\xi\xi_1(\xi-\xi_1)|\geq \frac 12|\sigma_1|$$ and
$$|\xi|\leq |\sigma_1|,$$ which combined with (\ref{calc1}) and smoothing relation (\ref{smooth}) yield
\begin{eqnarray*}
K_{422} &\lesssim & \frac
1{\langle\sigma_1\rangle^{1/2}}\Big(\int_{\widetilde{\Omega}_{412}}\frac{\langle\xi\xi_1(\xi-\xi_1)\rangle^{-2s}\langle\xi\rangle^{2+4s}}
{\langle\sigma\rangle^{1-\delta}\langle\sigma_2\rangle^{1-\delta}}d\tau
d\xi\Big)^{1/2}\\ &\lesssim &
\langle\sigma_1\rangle^{-s-1/2}\Big(\int_{\Omega_{412}'}
\frac{\langle\xi\rangle^{2+4s}}{\langle\sigma_1+3\xi\xi_1(\xi-\xi_1)\rangle^{1-2\delta}}d\xi
\Big)^{1/2}\\ &\lesssim &
\langle\sigma_1\rangle^{-s-1+\delta}\Big(\int_{|\xi|\leq
|\sigma_1|}\frac{d\xi}{\langle\xi\rangle^{-4s-2}}\Big)^{1/2},
\end{eqnarray*}
(we have set $\Omega_{412}'=\{\xi : \exists
\tau,\tau_1,\xi_1\in\R,
(\tau,\tau_1,\xi,\xi_1)\in\Omega_{412}\}$).\\
Now from the assumptions (\ref{conds}) on $s$ we see that
$$K_{422}\lesssim\langle\sigma_1\rangle^{-s-1+\delta}\langle\sigma_1\rangle^{2s+3/2}\lesssim
1$$ if $\alpha\leq 1/2$ and
$$K_{422}\lesssim\langle\sigma_1\rangle^{-s-1+\delta}\lesssim 1$$
otherwise.

\vskip 0.5cm \underline{Estimate in $\Omega_{413}$ :} In this
domain, $|\xi_1|\sim|\xi-\xi_1|$ thus using (\ref{est1}) it
follows that
\begin{eqnarray*}
K_{413} &\lesssim &
\frac{\langle\xi_1\rangle^{-s-1/2+\delta}}{\langle\sigma_1\rangle^{\delta}}
\Big(\int_{\widetilde{\Omega}_{413}}\frac{\langle\xi\rangle^{2s+1+2\delta}\langle\xi-\xi_1\rangle^{-2s-1+2\delta}}
{\langle\sigma\rangle^{1-\delta}\langle\sigma_2\rangle^{1-\delta}}d\tau
d\xi\Big)^{1/2}\\ &\lesssim &
\frac{\langle\xi_1\rangle^{-2s-1+2\delta}}{\langle\sigma_1\rangle^{\delta}}
\Big(\int_{\Omega_{413}'}\frac{d\xi}{\langle\sigma_1+3\xi\xi_1(\xi-\xi_1)\rangle^{1-2\delta}}\Big)^{1/2}.
\end{eqnarray*}
Following the works of Kenig, Ponce and Vega \cite{MR1329387}, we
perform the change of variables
$\mu_1=\sigma_1+3\xi\xi_1(\xi-\xi_1)$. Thus, since $\displaystyle
d\xi\sim\frac{d\mu_1}{|\xi_1|^{1/2}\sqrt{|4\tau_1-\xi_1^3-4\mu_1|}}$
and in view of (\ref{calc2}), we bound $K_{413}$ by
\begin{eqnarray*}
K_{413} &\lesssim &
\frac{\langle\xi_1\rangle^{-2s-5/4+2\delta}}{\langle\sigma_1\rangle^{\delta}}\Big(\int_{|\mu_1|\leq
2|\sigma_1|}\frac{d\mu_1}{\langle\mu_1\rangle^{1-2\delta}\sqrt{|4\tau_1-\xi_1^3-4\mu_1|}}\Big)^{1/2}\\
&\lesssim &
\frac{\langle\xi_1\rangle^{-2s-5/4+2\delta}}{\langle\sigma_1\rangle^{\delta}}\frac{\langle\sigma_1\rangle^{\delta}}{\langle
4\tau_1-\xi_1^3\rangle^{1/4}}\\ &\lesssim &
\langle\xi_1\rangle^{-2s-5/4+2\delta}\langle
4\tau_1-\xi_1^3\rangle^{-1/4}.
\end{eqnarray*}
Note that in $\Omega_{413}$, we have $|\sigma_1|\leq
\frac{12}{100}|\xi_1|^3$, which leads to
$$3|\xi_1|^3\leq|4\sigma_1+3\xi_1^3|+4|\sigma_1|\leq
|4\tau_1-\xi_1^3|+\frac{48}{100}|\xi_1|^3$$ and thus
$|\xi_1|^3\lesssim |4\tau_1-\xi_1^3|$. One deduce that for
$-2s-5/4+2\delta>0$,
$$K_{413}\lesssim \langle 4\tau_1-\xi_1^3\rangle^{\frac
13(-2s-5/4+2\delta)-1/4}\lesssim \langle
4\tau_1-\xi_1^3\rangle^{-\frac 23(s+1)+2\delta/3}\lesssim 1.
$$

\vskip 0.5cm \textbf{Estimate in $\Omega_{42}$}\\ We split this
region in two components :
\begin{eqnarray*} \Omega_{421} &=&
\{(\tau,\tau_1,\xi,\xi_1)\in\Omega_{42} : |\xi_1|\leq 100|\xi|\},\\
\Omega_{422} &=& \{(\tau,\tau_1,\xi,\xi_1)\in \Omega_{42} :
|\xi_1|\geq 100|\xi|\}.
\end{eqnarray*}

\vskip 0.5cm \underline{Estimate in $\Omega_{421}$ :} In
$\Omega_{421}$, $\langle\xi_1\rangle\lesssim \langle\xi\rangle$
and $\langle\xi-\xi_1\rangle\lesssim \langle\xi\rangle$. Then, we
bound $I_{421}$ exactly in the same way that for $I_{411}$. After
Plancherel and H\"{o}lder, we are lead to the estimate
\begin{eqnarray*}I_{421}  &\lesssim &
\|h\|_{L^2_{xt}}\Big\|\mathcal{F}^{-1}\Big(\frac{\langle\xi\rangle^{[-s/2-1/4-\alpha(1-\lambda)+3\delta/2]^+}\widehat{f}}
{\langle\sigma\rangle^{\lambda/2}}\Big)\Big\|_{L^4_{xt}}\\ &&
\times
\Big\|\mathcal{F}^{-1}\Big(\frac{\langle\xi\rangle^{[-s/2-1/4-\alpha(1-\lambda)+3\delta/2]^+}\widehat{g}}
{\langle\sigma\rangle^{\lambda/2}}\Big)\Big\|_{L^4_{xt}}.
\end{eqnarray*}
It suffices to choose $\lambda=3/4+\eps$ to apply Lemma \ref{stri}
with $s>-3/4-\alpha/2$.

\vskip 0.5cm \underline{Estimate in $\Omega_{422}$ :} We split
$\Omega_{422}$ into three sub-domains
\begin{eqnarray*} \Omega_{4221} &= &
\{(\tau,\tau_1,\xi,\xi_1)\in\Omega_{422} :
3|\xi\xi_1(\xi-\xi_1)|\leq \frac 12|\sigma|\},\\ \Omega_{4222} &=&
\{(\tau,\tau_1,\xi,\xi_1)\in\Omega_{422} :
3|\xi\xi_1(\xi-\xi_1)|\geq \frac 12|\sigma|, |\sigma_2|\leq 1\},\\
\Omega_{4223} &=& \{(\tau,\tau_1,\xi,\xi_1)\in\Omega_{422} :
3|\xi\xi_1(\xi-\xi_1)|\geq \frac 12|\sigma|, |\sigma_2|\geq 1\}.
\end{eqnarray*}

\vskip 0.5cm \underline{Estimate in $\Omega_{4221}$ :} In this
region one has that
$$|\sigma+3\xi\xi_1(\xi-\xi_1)|\geq \frac 12|\sigma|$$ and since $\sigma_1+\sigma_2-\sigma=3\xi\xi_1(\xi-\xi_1)$
it follows that
$$|\sigma_1|\leq |\sigma|\leq 2|\sigma_1+\sigma_2|\leq
4|\sigma_1|$$ and $|\sigma_1|\sim|\sigma|\geq
|\xi\xi_1(\xi-\xi_1|$. Therefore using (\ref{est1}), one obtain
\begin{eqnarray*}K_{4221} &\lesssim &
\frac{\langle\xi_1\rangle^{-s}}{\langle\sigma_1\rangle^{1/2}}\Big(
\int_{\widetilde{\Omega}_{4221}}\frac{|\xi|^2\langle\xi\rangle^{2s}\langle\xi-\xi_1\rangle^{-2s}}
{\langle\sigma\rangle^{1-\delta}\langle\sigma_2\rangle}d\tau
d\xi\Big)^{1/2}\\ &\lesssim &
\frac{1}{\langle\sigma_1\rangle^{1-\delta}}\Big(
\int_{\widetilde{\Omega}_{4221}}\frac{\langle\xi\xi_1(\xi-\xi_1)\rangle^{-2s}\langle\xi\rangle^{2+4s}}
{\langle\sigma_2\rangle^{1+\delta}}d\tau d\xi\Big)^{1/2}\\
&\lesssim &
\langle\sigma_1\rangle^{-s-1+\delta}\Big(\int_{|\xi|\lesssim |\sigma_1|}\frac{d\xi}{\langle\xi\rangle^{-4s-2}}\Big)^{1/2}\\
& \lesssim & 1
\end{eqnarray*}
as for $K_{422}$.

\vskip 0.5cm \underline{Estimate in $\Omega_{4222}$ :} First
consider the case $\alpha\leq 1/2$. Then,
\begin{eqnarray*} K_{4222} &\lesssim &
\frac{|\xi|\langle\xi\rangle^s}{\langle\sigma\rangle^{1/2-\delta/2}}\Big(\int_{\widetilde{\Omega}_{4222}}\frac{\langle\xi_1\rangle^{-2s}
\langle\xi-\xi_1\rangle^{-2s}}{\langle\sigma_1\rangle\langle\sigma_2\rangle}d\tau_1
d\xi_1\Big)^{1/2}\\ &\lesssim &
\frac{|\xi|^{1+s}\langle\xi\rangle^s}{\langle\sigma\rangle^{1/2-\delta/2}}\Big(\int_{\widetilde{\Omega}_{4222}}\frac
{\langle\xi\xi_1(\xi-\xi_1)\rangle^{-2s}}{\langle\sigma_1\rangle^{1-\delta}\langle\sigma_2\rangle^{1-\delta}}d\tau_1d\xi_1\Big)^{1/2}\\
&\lesssim &
|\xi|^{1+s}\langle\xi\rangle^s\langle\sigma\rangle^{-s-1/2+\delta/2}\Big(\int_{\Omega_{4222}'}\frac{d\xi_1}
{\langle\sigma+3\xi\xi_1(\xi-\xi_1)\rangle^{1-2\delta}}
\Big)^{1/2}.
\end{eqnarray*}
The change of variables $\mu=\sigma+3\xi\xi_1(\xi-\xi_1)$ gives
the inequalities
\begin{eqnarray*}K_{4222} &\lesssim &
|\xi|^{3/4+s}\langle\xi\rangle^s\langle\sigma\rangle^{-s-1/2+\delta/2}\Big(\int_{|\mu|\leq
2|\sigma|}\frac{d\mu}{\langle\mu\rangle^{1-2\delta}\sqrt{|4\tau-\xi^3-4\mu|}}\Big)^{1/2}\\
&\lesssim &
|\xi|^{3/4+s}\langle\xi\rangle^s\langle\sigma\rangle^{-s-1/2+3\delta/2}\langle4\tau-\xi^3\rangle^{-1/4},
\end{eqnarray*}
which is bounded on $\R^2$. If  $\alpha>1/2$, we have directly
\begin{eqnarray*} K_{4222} &\lesssim &
\frac{|\xi|\langle\xi\rangle^s}{\langle\sigma\rangle^{1/2-\delta/2}}
\Big(\int_{\widetilde{\Omega}_{4222}}
\frac{\langle\xi_1\rangle^{-4s}}{\langle\xi_1\rangle^{2\alpha}\langle\xi-\xi_1\rangle^{2\alpha}\langle\sigma_2\rangle^{1+\eps}}
d\tau_1d\xi_1\Big)^{1/2}\\ &\lesssim &
\langle\xi\rangle^{1/2+s+\delta/2}\Big(\int_\R\frac{d\xi_1}{\langle\xi_1\rangle^{4s+2+4\alpha+2\delta}}\Big)^{1/2}\\
&\lesssim & 1
\end{eqnarray*}
for $s>-1/4-\alpha$.

\vskip 0.5cm \underline{Estimate in $\Omega_{4223}$ :} As we will
see below, that is in this sub-domain that the condition
$s>s_\alpha$ appears. Also, to obtain our estimates, we will need
tu use a dyadic decomposition of the variables $\xi_1$,
$\xi-\xi_1$, $\xi$ and $\sigma_1$, $\sigma_2$, $\sigma$.  Hence,
following the notations introduced in Subsection \ref{subsec-tao},
we have to bound

\begin{eqnarray*}\lefteqn{I_{4223} = \sum_{\substack{N_3\ll N_1\sim
N_2\\ N_1\gtrsim 1}}\sum_{\substack{L_3\gtrsim L_1\gtrsim L_2\\
L_3\sim N_3N_1^2}}\frac{N_3\langle
N_3\rangle^sN_1^{-s}N_2^{-s}}{\max(L_3,N_3^{2\alpha})^{1/2-2\delta}\max(L_1,N_1^{2\alpha})^{1/2}\max(L_2,N_2^{2\alpha})^{1/2}}}\\
& \times \displaystyle
\int_{\R^4}\frac{\widehat{h}(\tau,\xi)}{\langle\sigma\rangle^\delta}\widehat{f}(\tau_1,\xi_1)\widehat{g}(\tau-\tau_1,\xi-\xi_1)\chi_{|\xi|\sim
N_3, |\xi_1|\sim N_1, |\xi-\xi_1|\sim N_2}\chi_{|\sigma|\sim L_3,
|\sigma_1|\sim L_1, |\sigma_2|\sim L_2} d\tau d\tau_1 d\xi d\xi_1.
\end{eqnarray*}
Using the (+-) coherence case of Lemma \ref{tao} as well as Lemma
\ref{contract}, we get
\begin{eqnarray*}\lefteqn{\int_{\R^4}\frac{\widehat{h}(\tau,\xi)}{\langle\sigma\rangle^\delta}
\widehat{f}(\tau_1,\xi_1)\widehat{g}(\tau-\tau_1,\xi-\xi_1)\chi_{|\xi|\sim
N_3, |\xi_1|\sim N_1, |\xi-\xi_1|\sim N_2}\chi_{|\sigma|\sim L_3,
|\sigma_1|\sim L_1, |\sigma_2|\sim L_2} d\tau d\tau_1 d\xi
d\xi_1}\\  && \lesssim
L_2^{1/2}N_1^{-1}\min\big(N_3N_1^2,\frac{N_1}{N_3}L_1\big)^{1/2}
\Big\|\mathcal{F}^{-1}\Big(\frac{\widehat{h}(\tau,\xi)}{\langle\sigma\rangle^\delta}\Big)\Big\|_{L^2_{xt}}\|f\|_{L^2_{xt}}\|g\|_{L^2_{xt}}\\
 && \lesssim T^\nu
L_2^{1/2}N_1^{-1}\min\big(N_3N_1^2,\frac{N_1}{N_3}L_1\big)^{1/2}
\|h\|_{L^2_{xt}}\|f\|_{L^2_{xt}}\|g\|_{L^2_{xt}}.
\end{eqnarray*}
Thus we reduce to show \begin{equation}\label{f-sum}\begin{split}
\sum_{\substack{N_3\ll N_1\sim
N_2\\ N_1\gtrsim 1}}\sum_{\substack{L_3\gtrsim L_1\gtrsim L_2\\
L_3\sim N_3N_1^2}} \frac{N_3\langle
N_3\rangle^sN_1^{-s}N_2^{-s}}{\max(L_3,N_3^{2\alpha})^{1/2-2\delta}\max(L_1,N_1^{2\alpha})^{1/2}\max(L_2,N_2^{2\alpha})^{1/2}}\\
 \times
L_2^{1/2}N_1^{-1}\min\big(N_3N_1^2,\frac{N_1}{N_3}L_1\big)^{1/2}\lesssim
1.
\end{split}
\end{equation}
Recalling that $N_j=2^{n_j}$ and $L_j=2^{\ell_j}$, $j=1,2,3$, for
all $\lambda\in[0,1]$, the right hand side of (\ref{f-sum}) is
bounded by

\begin{eqnarray*} &\lesssim &\sum_{\substack{N_3\ll N_1\sim
N_2\\ N_1\gtrsim 1}}\sum_{\substack{L_3\gtrsim L_1\gtrsim L_2\\
L_3\sim N_3N_1^2}}\frac{N_3\langle N_3\rangle^s N_1^{-2s-1}
(N_3N_1^2)^{\lambda/2}\big(\frac{N_1}{N_3}L_1\big)^{(1-\lambda)/2}}{L_3^\delta
(N_3N_1^2)^{1/2-3\delta}L_1^{(1-\lambda)/2}N_1^{\alpha\lambda}}\\
&\lesssim & \sum_{\substack{N_3\ll N_1\sim N_2\\ N_1\gtrsim 1}}
\Big(\sum_{L_1,L_2,L_3\gtrsim
1}\frac{1}{(L_1L_2L_3)^{\delta/3}}\Big)N_3^{\lambda+3\delta}\langle
N_3\rangle^s N_1^{-2s-3/2+\lambda(1/2-\alpha)+6\delta}\\ &\lesssim
& \sum_{\substack{N_1\gtrsim 1\\ N_3\ll
N_1}}N_3^{\lambda+3\delta}\langle N_3\rangle^s
N_1^{-2s-3/2+\lambda(1/2-\alpha)+6\delta}.
\end{eqnarray*}
This last expression is finite if $\lambda+3\delta<-s$ and
$-2s-3/2+\lambda(1/2-\alpha)+6\delta<0$.\\
When $\alpha\leq 1/2$, it suffices to choose $\lambda=0$ and
$\lambda=-s-4\delta$ otherwise.

\section{Proof of the main result}
\label{sec-main} In this section, we briefly indicate how the
results stated in Section \ref{subsec-lin} and the bilinear
estimate (\ref{bil}) yield Theorem \ref{main} (see for instance
\cite{MR1918236} for the details).

Actually, local existence of a solution is a consequence of the
following modified version of Theorem \ref{thbil}.

\begin{proposition} Given $s_c^+>s_\alpha$, there exist
$\nu,\delta>0$ such that for any $s\geq s_c^+$ and any $u,v\in
X^{1/2,s}_\alpha$ with compact support in $[-T,+T]$,
\begin{equation}\label{modbil}\|\partial_x(uv)\|_{X^{-1/2+\delta,s}_\alpha}\lesssim
T^\nu(\|u\|_{X^{1/2,s_c^+}_\alpha}\|v\|_{X^{1/2,s}_\alpha}+\|u\|_{X^{1/2,s}_\alpha}\|v\|_{X^{1/2,s_c^+}_\alpha}).\end{equation}
\end{proposition}
Estimate (\ref{modbil}) is obtained thanks to (\ref{bil}) and the
triangle inequality $$\forall s\geq s_c^+,\
\langle\xi\rangle^s\leq
\langle\xi\rangle^{s_c^+}\langle\xi_1\rangle^{s-s_c^+}+\langle\xi\rangle^{s_c^+}\langle\xi-\xi_1\rangle^{s-s_c^+}.$$
Let $u_0\in H^s(\R)$ with $s$ obeying (\ref{conds}). Define $F(u)$
as
$$F(u)=F_{u_0}(u)=\psi(t)\Big[W_\alpha(t)u_0-\frac{\chi_{\R_+}(t)}{2}\int_0^tW_\alpha(t-t')\partial_x(\psi_T^2(t')
u^2(t'))dt'\Big].$$ We shall prove that for $T\ll 1$, $F$ is
contraction in a ball of the Banach space $$Z=\{u\in
X^{1/2,s}_\alpha,\
\|u\|_Z=\|u\|_{X^{1/2,s_c^+}_\alpha}+\gamma\|u\|_{X^{1/2,s}_\alpha}<+\infty\},$$
where $\gamma$ is defined for all nontrivial $\varphi$ by
$$\gamma=\frac{\|\varphi\|_{H^{s_c^+}}}{\|\varphi\|_{H^s}}.$$
Combining (\ref{lin-free}), (\ref{lin-for}) as well as
(\ref{modbil}), it is easy to derive that $$\|F(u)\|_Z\leq
C(\|u_0\|_{H^{s_c^+}}+\gamma\|u_0\|_{H^s})+CT^\nu\|u\|_Z^2$$ and
$$\|F(u)-F(v)\|_Z\leq CT^\nu\|u-v\|_Z\|u+v\|_Z$$ for some
$C,\nu>0$. Thus, taking $T=T(\|u_0\|_{H^{s_c^+}})$ small enough,
we deduce that $F$ is contractive on the ball of radius
$4C\|u_0\|_{H^{s_c^+}}$ in $Z$. This proves the existence of a
solution $u$ to $u=F(u)$ in $X^{1/2,s}_{\alpha,T}$.

Following similar arguments of \cite{MR1918236}, it is not too
difficult to see that if $u_1,u_2\in X^{1/2,s}_{\alpha, T}$ are
solutions of (\ref{int}) and $0<\delta<T/2$, then there exists
$\nu>0$ such that
$$\|u_1-u_2\|_{X^{1/2,s}_{\alpha,\delta}}\lesssim
T^\nu\big(\|u_1\|_{X^{1/2,s}_{\alpha,T}}+\|u_2\|_{X^{1/2,s}_{\alpha,T}}\big)\|u_1-u_2\|_{X^{1/2,s}_{\alpha,\delta}},$$
which leads to $u_1\equiv u_2$ on $[0,\delta]$, and then on
$[0,T]$ by iteration. This proves the uniqueness of the solution.

It is straightforward to check that
$W_\alpha(\cdot)u_0\in\C(\R_+,H^s(\R))\cap\C(\R_+^\ast,H^\infty(\R))$.
Then it follows from Theorem \ref{thbil}, Lemma \ref{glob} and the
local existence of the solution that $$u\in\C([0,T],
H^s(\R))\cap\C((0,T],H^\infty(\R))$$ for some
$T=T(\|u_0\|_{H^{s_c^+}})$. By induction, we have
$u\in\C((0,T],H^\infty(\R))$. Taking the $L^2$-scalar product of
(\ref{eq}) with $u$, we obtain that $t\mapsto\|u(t)\|_{H^{s_c^+}}$
is nonincreasing on $(0,T]$. Since the existence time of the
solution depends only on the norm $\|u_0\|_{H^{s_c^+}}$, this
implies that the solution can be extended globally in time.

\section{Proof of Theorem \ref{opt}}
\label{sec-opt} Let $s\geq s_\alpha$. To prove Theorem \ref{opt},
it suffices to show that the multiplier
$$\mathcal{M}=\sup_{\substack{f,g,h\in\S(\R^2)\\\|f\|_{L^2_{xt}},\|g\|_{L^2_{xt}},\|h\|_{L^2_{xt}}}\lesssim
1}\Big|\int_{\R^4}K(\tau,\tau_1,\xi,\xi_1)\widehat{h}(\tau,\xi)\widehat{f}(\tau_1,\xi_1)\widehat{g}(\tau-\tau_1,\xi-\xi_1)
d\tau d\tau_1d\xi d\xi_1\Big|$$ is infinite. Setting
$$\mathcal{A}=\{f\in \S(\R^2):\ \widehat{f}\geq 0,\
 \|f\|_{L^2_{xt}}\lesssim 1\},$$
 and performing as previously a dyadic decomposition of the variables $\xi_1$, $\xi-\xi_1$, $\xi$
and $\sigma_1$, $\sigma_2$, $\sigma$, it follows that
\begin{eqnarray*}\mathcal{M} &\gtrsim & \sup_{f,g,h\in\mathcal{A}} \sum_{N_1,N_2,N_3}\sum_{L_1,L_2,L_3}\frac{N_3\langle
N_3\rangle^s\langle N_1\rangle^{-s}\langle
N_2\rangle^{-s}}{\langle\max(L_3,N_3^{2\alpha})\rangle^{1/2-\delta}\langle\max(L_1,N_1^{2\alpha})\rangle^{1/2}\langle
\max(L_2,N_2^{2\alpha})\rangle^{1/2}}\\ && \times
\int_{\R^4}\widehat{h}(\tau,\xi)
\widehat{f}(\tau_1,\xi_1)\widehat{g}(\tau-\tau_1,\xi-\xi_1)\chi_{|\xi|\sim
N_3, |\xi_1|\sim N_1, |\xi-\xi_1|\sim N_2}\chi_{|\sigma|\sim L_3,
|\sigma_1|\sim L_1, |\sigma_2|\sim L_2} \\ &\gtrsim &
\sup_{N_1,N_2,N_3}\sup_{L_1,L_2,L_3}\frac{N_3\langle
N_3\rangle^s\langle N_1\rangle^{-s}\langle
N_2\rangle^{-s}}{\langle\max(L_3,N_3^{2\alpha})\rangle^{1/2-\delta}\langle\max(L_1,N_1^{2\alpha})\rangle^{1/2}\langle
\max(L_2,N_2^{2\alpha})\rangle^{1/2}}\\ && \times
\sup_{f,g,h\in\mathcal{A}} \int_{\R^4}\widehat{h}(\tau,\xi)
\widehat{f}(\tau_1,\xi_1)\widehat{g}(\tau-\tau_1,\xi-\xi_1)\chi_{|\xi|\sim
N_3, |\xi_1|\sim N_1, |\xi-\xi_1|\sim N_2}\chi_{|\sigma|\sim L_3,
|\sigma_1|\sim L_1, |\sigma_2|\sim L_2}.
\end{eqnarray*}
Now we localize the previous supremum to the critical region
$$\left\{\begin{array}{ll}N_3\ll N_1\sim N_2\\
N_3N_1^2\sim L_3\gtrsim L_1\gtrsim L_2\gtrsim 1
\end{array}\right.$$ which corresponds to a sub-domain of
$\Omega_{4223}$. In this case, the optimality of (\ref{+-}) gives
\begin{eqnarray*}\sup_{f,g,h\in\mathcal{A}}
\int_{\R^4}\widehat{h}(\tau,\xi)
\widehat{f}(\tau_1,\xi_1)\widehat{g}(\tau-\tau_1,\xi-\xi_1)\chi_{|\xi|\sim
N_3, |\xi_1|\sim N_1, |\xi-\xi_1|\sim N_2}\chi_{|\sigma|\sim L_3,
|\sigma_1|\sim L_1, |\sigma_2|\sim L_2}\\ \gtrsim
L_2^{1/2}N_1^{-1}\min\big(N_3N_1^2,\frac{N_1}{N_3}L_1\big)^{1/2}.
\end{eqnarray*}
Therefore we have the bound
\begin{eqnarray*} \mathcal{M} &\gtrsim & \sup_{1\lesssim N_3\ll
N_1\sim N_2}\sup_{\substack{1\lesssim L_2\lesssim L_1\lesssim
L_3\\ L_3\sim N_3N_1^2}}\frac{N_3^{1+s}N_1^{-2s}
L_2^{1/2}N_1^{-1}\min\big(N_3N_1^2,\frac{N_1}{N_3}L_1\big)^{1/2}}{\max(N_3N_1^2,N_3^{2\alpha})^{1/2-\delta}\max(L_1,N_1^{2\alpha})^{1/2}
\max(L_2,N_2^{2\alpha})^{1/2}}.
\end{eqnarray*}
First consider the case $0\leq \alpha\leq 1/2$. Then,
\begin{eqnarray*}\mathcal{M} &\gtrsim & \sup_{\substack{N_1\gg 1\\
N_3\sim 1}}\sup_{L_1\sim L_2\sim
N_3^2N_1}\frac{N_1^{-2s-1}L_2^{1/2}\big(\frac{N_1}{N_3}L_1\big)^{1/2}}{N_1^{1-2\delta}L_1^{1/2}L_2^{1/2}}\\
&\gtrsim & \sup_{N_1\gg 1} N_1^{-2s-3/2+2\delta}\\ &=&+\infty
\end{eqnarray*}
for $s\leq -3/4$. Now if $1/2<\alpha\leq 1$, we estimate
\begin{eqnarray*}\mathcal{M} &\gtrsim & \sup_{\substack{N_1\gg 1\\
N_3\sim N_1^{\alpha-1/2}}}\sup_{L_1\sim L_2\sim
N_1^{2\alpha}}\frac{N_3^{1+s}N_1^{-2s-1}L_2^{1/2}N_1^{3/4+\alpha/2}}{(N_3N_1^2)^{1/2-\delta}L_2^{1/2}N_1^\alpha}\\
&\gtrsim & \sup_{N_1\gg
1}N_1^{-2s-5/4-\alpha/2+(\alpha-1/2)(s+1/2)+(\alpha+3/2)\delta}\\
&=& +\infty
\end{eqnarray*}
for $-2s-5/4-\alpha/2+(\alpha-1/2)(s+1/2)\geq 0$, i.e.
$\displaystyle s\leq \frac{-3}{5-2\alpha}$.

\bibliographystyle{plain}
\bibliography{ref}

\end{document}